\documentclass[12pt]{amsart}
\usepackage{amssymb,latexsym,color}
\usepackage{xspace}
\numberwithin{equation}{section}

\newtheorem{theorem}{Theorem}[section]
\newtheorem{proposition}[theorem]{Proposition}

\newtheorem{corollary}[theorem]{Corollary}
\newtheorem{lemma}[theorem]{Lemma}

\theoremstyle{definition}

\newtheorem{remark}[theorem]{Remark}
\newtheorem{example}[theorem]{Example}
\newtheorem{definition}[theorem]{Definition}

\renewcommand{\eqref}[1]{{\rm (\ref{#1})}}

\def\proof{\smallskip\noindent {\bf Proof.\ }}
\def\endproof{\hfill$\square$\medskip}

\def\AA{\mathbb{A}}
\def\BB{\mathbb{B}}
\def\CC{\mathbb{C}}
\def\DD{\mathbb{D}}
\def\EE{\mathbb{E}}
\def\FF{\mathbb{F}}
\def\GG{\mathbb{G}}

\def\RR{\mathbb{R}}
\def\ZZ{\mathbb{Z}}
\def\sgn{\operatorname{sgn}}

\def\Gam{\Gamma}

\def\qCm{quasi-Cartan matrix\xspace}
\def\qCms{quasi-Cartan matrices\xspace}

\def\CO{cyclically orientable\xspace}
\def\COP{cyclical orientability\xspace}
\def\COD{cyclically oriented\xspace}
\def\QF{quasi-finite\xspace}
\def\QC{quasi-Cartan\xspace}

\newcommand{\Con}{\operatorname{Con}}
\newcommand{\Ver}{\operatorname{Ver}}
\newcommand{\Edg}{\operatorname{Edg}}

\newcommand{\Hol}{\operatorname{Cyc}}
\newcommand{\holes}{chordless cycles\xspace}
\newcommand{\hole}{chordless cycle\xspace}

\newcommand{\HVCenter}[1]{\setbox 0=\hbox{#1}%
        \dimen0=\wd0%
        \dimen1=\ht0%
        \divide\dimen0 by 2%
        \divide\dimen1 by 2%
        \hskip -\dimen0%
        \lower \dimen1%
        \box0%
        \hskip -\dimen0}
\newcommand{\HBCenter}[1]{\setbox 0=\hbox{#1}%
        \dimen0=\wd0%
        \dimen1=\ht0%
        \divide\dimen0 by 2%
        \hskip -\dimen0%
        \box0%
        \hskip -\dimen0}
\newcommand{\HTCenter}[1]{\setbox 0=\hbox{#1}%
        \dimen0=\wd0%
        \dimen1=\ht0%
        \divide\dimen0 by 2%
        \hskip -\dimen0%
        \lower \dimen1%
        \box0%
        \hskip -\dimen0}
\newcommand{\RVCenter}[1]{\setbox 0=\hbox{#1}%
        \dimen0=\wd0%
        \dimen1=\ht0%
        \divide\dimen1 by 2%
        \hskip -\dimen0%
        \lower \dimen1%
        \box0%
        \hskip -\dimen0}
\newcommand{\RTCenter}[1]{\setbox 0=\hbox{#1}%
        \dimen0=\wd0%
        \dimen1=\ht0%
        \hskip -\dimen0%
        \lower \dimen1%
        \box0%
        \hskip -\dimen0}
\newcommand{\LTCenter}[1]{\setbox 0=\hbox{#1}%
        \dimen1=\ht0%
        \lower \dimen1%
        \box0%
        \hskip -\dimen0}
\newcommand{\RBCenter}[1]{\setbox 0=\hbox{#1}%
        \dimen0=\wd0%
        \hskip -\dimen0%
        \box0%
        \hskip -\dimen0}
\newcommand{\LVCenter}[1]{\setbox 0=\hbox{#1}%
        \dimen1=\ht0%
        \divide\dimen1 by 2%
        \lower \dimen1%
        \box0%
        \hskip -\dimen0}
\begin{document}

\title[Cluster algebras and positive matrices]
{Cluster algebras of finite type and positive symmetrizable
matrices}
\date{December~13, 2004; revised June 27, 2005.}

\author{Michael Barot}
\address{\noindent Instituto de Matem\'aticas,  UNAM,
Ciudad Universitaria, 04510 Mexico D.F., Mexico}
\email{barot@matem.unam.mx}

\author{Christof Geiss}
\address{\noindent Instituto de Matem\'aticas,  UNAM,
Ciudad Universitaria, 04510 Mexico D.F., Mexico}
\email{christof@matem.unam.mx}

\author{Andrei Zelevinsky}
\address{\noindent Department of Mathematics, Northeastern University,
 Boston, MA 02115, USA}
\email{andrei@neu.edu}

\begin{abstract}
The paper is motivated by an analogy between cluster algebras and Kac-Moody
algebras: both theories share the same classification of  finite type
objects by familiar Cartan-Killing types. However the underlying
combinatorics beyond the two classifications is different: roughly
speaking, Kac-Moody algebras are associated with  (symmetrizable)
Cartan matrices, while cluster algebras correspond to
skew-symmetrizable matrices. We study an interplay between the two
classes of matrices, in particular, establishing a new criterion
for deciding whether a given skew-symmetrizable matrix
gives rise to a cluster algebra of finite type.
\end{abstract}

 \thanks{Research
supported
by DGAPA grant IN101402-3 (C.G.) and NSF (DMS) grant 0200299~(A.Z.).}

\maketitle

\section{Introduction}
This paper is motivated by the theory of cluster algebras
initiated in~\cite{FZ1}.
Here we deal exclusively with the combinatorial aspects of the theory,
so no knowledge of algebraic properties of cluster algebras
(including their definition) will be assumed or needed.
The reader should just bear in mind an analogy between cluster algebras
and Kac-Moody algebras.
In both theories, there is an appropriate notion of finite type.
(For Kac-Moody algebras, ``finite type" just means being
finite-dimensional, that is, a semisimple Lie algebra.)
Cluster algebras of finite type were classified in~\cite{FZ2},
and the resulting classification turns out to be identical to the famous
Cartan-Killing classification of semisimple Lie algebras.
However the underlying combinatorics beyond the two classes of
algebras is different: roughly speaking, Kac-Moody algebras
correspond to (symmetrizable) Cartan matrices,
while cluster algebras correspond to skew-symmetrizable matrices.
In this paper, we study an interplay between the two classes of matrices.
In particular, we establish a new criterion for deciding
whether a given skew-symmetrizable matrix gives rise to a cluster
algebra of finite type.

To state our main results, we need some terminology.
In what follows, by a matrix we always mean a square integer matrix.
A matrix~$A$ (resp.~$B$) is symmetrizable (resp.
skew-symmetrizable) if $DA$ (resp.~$DB$) is symmetric (resp.
skew-symmetric) for some diagonal matrix~$D$ with positive
diagonal entries.
Thus, in a symmetrizable (resp. skew-symmetrizable) matrix, the
two transpose entries have the same sign (resp. opposite signs).
We say that a symmetrizable matrix is a \emph{\qCm} if all its diagonal entries
are equal to~$2$ (note that a skew-symmetrizable matrix
has all diagonal entries equal to~$0$).
We say that a \qCm~$A$ is \emph{positive} if the
symmetrized matrix~$DA$ is positive definite; by the Sylvester
criterion, this means that the principal minors of~$A$ are all positive.
For a skew-symmetrizable matrix~$B$,
we will refer to a \qCm~$A$ with $|A_{ij}| = |B_{ij}|$
for all~$i \neq j$ as a \emph{\QC companion} of~$B$.

A \qCm is a (generalized) \emph{Cartan matrix} if its off-diagonal entries
are non-positive. These are the matrices giving rise to Kac-Moody algebras,
see~\cite{kac}. As shown in \cite{kac}, the Kac-Moody algebra associated
to~$A$ is finite-dimensional if and only if~$A$ is positive.
Thus, positive Cartan matrices form the combinatorial backbone of the
Cartan-Killing classification: every such matrix can be
transformed by a simultaneous permutation of rows and columns to
a block-diagonal matrix whose diagonal blocks are of the familiar
types~$\AA_n, \BB_n, \CC_n, \DD_n, \EE_6, \EE_7, \EE_8, \FF_4$, and~$\GG_2$,
represented by Dynkin diagrams.

On the other hand, cluster algebras are associated with the
mutation-equivalence classes of skew-symmetrizable matrices.
Recall  from~\cite{FZ1} that, for each matrix index~$k$,
the \emph{mutation} in direction~$k$ transforms a
skew-symmetrizable matrix~$B$ into another skew-symmetrizable
matrix $B' = \mu_k(B)$, whose entries are given by
\begin{equation}
\label{eq:matrix-mutation}
B'_{ij} =
\begin{cases}
-B_{ij} & \text{if $i=k$ or $j=k$;} \\
B_{ij} + \sgn(B_{ik}) [B_{ik} B_{kj}]_+
 & \text{otherwise,}
\end{cases}
\end{equation}
where we use the notation $[x]_+ = \max(x,0)$ and
$\sgn(x) = x/|x|$, with the convention $\sgn(0) = 0$
(the formula \eqref{eq:matrix-mutation} is easily seen to be equivalent
to~\cite[(4.3)]{FZ1}).
One easily checks that~$\mu_k$ is involutive, implying that
the repeated mutations in all directions give rise to the
\emph{mutation-equivalence} relation
on skew-symmetrizable matrices.
We are now ready to state the classification result from~\cite{FZ2}.

\begin{theorem}
\label{thm:FZ}
For a mutation-equivalence class~$\mathcal{S}$ of skew-symmetrizable matrices,
the following are equivalent.
\begin{itemize}
\item[\rm (1)]
The cluster algebra associated to~$\mathcal{S}$ is of finite type.
\item[\rm (2)]
$\mathcal{S}$ contains a matrix~$B$ such that the Cartan matrix~$A$
with off-diagonal entries~$A_{ij} = -|B_{ij}|$ is positive.
\item[\rm (3)]
For every~$B \in \mathcal{S}$ and all $i\neq j$, we have
$|B_{ij}B_{ji}| \leq 3$.
\end{itemize}
Furthermore, the Cartan-Killing type of the Cartan matrix~$A$ in (2) is
uniquely determined by~$\mathcal{S}$.
\end{theorem}

In this paper we address the following

\medskip

\noindent \emph{Recognition Problem.}
Given a skew-symmetrizable matrix~$B$, find an efficient way to
determine whether the cluster algebra associated with (the
mutation-equivalence class of)~$B$ is of finite type.

\medskip

The problem makes sense since the mutations are hard to control,
so each of the conditions (2) and (3) in Theorem~\ref{thm:FZ} is
hard to check in general.
A nice solution of the problem was obtained by
A.~Seven in~\cite{Seven}.
His answer is given in terms of ``forbidden minors" of~$B$.
Our answer is very different and probably more practical: it is
based on extending the criterion in~(2) to \emph{every}
representative of a mutation class in question.
To state it, we need a little bit more terminology.

To a skew-symmetrizable $n \times n$ matrix~$B$ we associate
a directed graph $\Gam(B)$ with vertices $1, \dots, n$ and
directed edges $i \to j$ for all $i,j$ with $B_{ij} > 0$.
A \emph{\hole} in $\Gam(B)$ is an induced subgraph isomorphic to a cycle
(thus, the vertices of a \hole can be labeled by the elements of
$\ZZ/p\ZZ$ for some $p \geq 3$ so that the edges between them
are precisely $\{i, i+1\}$ for $i \in \ZZ/p\ZZ$).

We are finally ready to state our main result.

\begin{theorem}
\label{thm:MainResult}
The mutation-equivalence class of a skew-symmetrizable matrix~$B$
satisfies the equivalent conditions in Theorem~\ref{thm:FZ} if and only if
$B$ satisfies:
\begin{itemize}
\item[\rm (4)]
Every \hole in $\Gam(B)$ is cyclically oriented, and~$B$
has a positive \QC companion.
\end{itemize}
\end{theorem}

\begin{remark}
\label{rmk:MainResult1}
We actually prove (2) $\Rightarrow$ (4)  $\Rightarrow$ (3), which
gives a new proof of the implication (2) $\Rightarrow$ (3), more
elementary and straightforward than the one in~\cite{FZ2}.
\end{remark}

Condition~(4) in Theorem~\ref{thm:MainResult} is not completely
explicit since~$B$ can have many \QC companions.
Note that a \QC companion of~$B$ is specified by choosing
the signs of its off-diagonal matrix entries,
with the only requirement that $\sgn(A_{ij}) = \sgn(A_{ji})$ for $i \neq j$.
Thus, the number of choices for~$A$ is $2^N$, where~$N$ is the
number of edges in $\Gam(B)$.
The following two propositions allow us to considerably sharpen
Theorem~\ref{thm:MainResult}.

\begin{proposition}
\label{pr:cyclic-A}
To be positive, a \QC companion~$A$ of
a skew-sym\-met\-rizable matrix~$B$ must satisfy the following sign condition:
for every \hole~$Z$ in $\Gam(B)$, the product
$\prod_{\{i,j\} \in Z} (-A_{ij})$
over all edges of~$Z$ must be negative.
\end{proposition}

\begin{proposition}
\label{pr:cyclic-B-to-A}
If every \hole in $\Gam(B)$ is cyclically oriented,
then~$B$ has a \QC companion (not necessarily positive) satisfying the sign
condition in Proposition~\ref{pr:cyclic-A}; furthermore, such a \QC companion
is unique up to simultaneous sign changes in rows and columns.
\end{proposition}

Proposition~\ref{pr:cyclic-A} (resp.~Proposition~\ref{pr:cyclic-B-to-A})
is a consequence of Proposition~\ref{lem:positive_cycles}
(resp. Corollary~\ref{cor:product-signs}) below.
In view of these results, in checking
condition~(4), it is enough to test positivity of just one
\QC companion of~$B$ since simultaneous sign changes in rows and columns
do not affect positivity.
The following example provides an illustration.

\begin{example}
Let $B(n)$ be the $n \times n$ skew-symmetric matrix with
the above-diagonal entries given by
\begin{equation}
\label{eq:B(n)}
B_{ij} =
\begin{cases}
-1 & \text{if $j-i=1$;} \\
1 & \text{if $j-i=2$;} \\
0 & \text{if $j-i > 2$.}
\end{cases}
\end{equation}

The graph $\Gam(B(n))$ has $n-2$ \holes: they are formed
by all triples of consecutive indices.
An immediate check shows that all of them are cyclically oriented.
Now let~$A(n)$ be the \QC companion of~$B(n)$ such that $A_{ij} =
B_{ij}$ for $j > i$.
An immediate check shows that $A(n)$ satisfies the condition in
Proposition~\ref{pr:cyclic-A}.
Let $d_n = \det(A(n))$, with the convention~$d_0 = 1$.
It is not hard to show that the generating function of this
sequence is given by
\begin{equation}
\label{eq:detAn}
\sum_{n \geq 0} d_n x^n = \frac{(1+x)(1+x+x^2)(1+x^2)(1+x^3)}
{1-x^{12}} \ ,
\end{equation}
implying that $d_{n+12} = d_n$ for $n \geq 0$.
Since the numerator in \eqref{eq:detAn} is a polynomial of degree~$8$,
we see that $d_9 = d_{10} = d_{11} = 0$.
The values of~$d_n$ for $1 \leq n \leq 8$ are given by the
following table:
\begin{table}[ht]
\begin{center}
\begin{tabular}{|c|c|c|c|c|c|c|c|c|}
\hline
&&&&&&&&\\[-.1in]
$n$ & 1 & 2 & 3 & 4 & 5 & 6 & 7  & 8\\
\hline
&&&&&&&&\\[-.1in]
$\det(A(n))$ &  2 & 3 & 4 & 4 & 4 & 3 & 2  & 1\\
\hline
&&&&&&&&\\[-.1in]
Cartan-Killing type &  $\AA_1$ & $\AA_2$ & $\AA_3$ & $\DD_4$ & $\DD_5$ &
$\EE_6$ & $\EE_7$  & $\EE_8$ \\[.05in]
\hline
\end{tabular}
\end{center}
\medskip
\caption{Determinants and Cartan-Killing types of the~$A(n)$}
\label{tab:det-An}
\end{table}

By the Sylvester criterion, $A(n)$ is positive if and only if $n \leq 8$.
Applying Theorem~\ref{thm:MainResult} and
Proposition~\ref{pr:cyclic-A}, we conclude that the cluster algebra
associated to $B(n)$ is of finite type precisely when $n \leq 8$.
The corresponding Cartan-Killing types are given in the last line
of Table~\ref{tab:det-An} (their determination is left as an exercise
to the reader).
\end{example}

The paper is organized as follows.
The proof of Theorem \ref{thm:MainResult} is carried out in the
next three sections.
In Section~\ref{sec:pos_properties}, which can be viewed as
a ``symmetrizable analogue'' of~\cite[Section~9]{FZ2}, we
establish some needed properties of positive \qCms.
In Section~\ref{sec:mutations}, we show that, under some conditions, the
mutation-equivalence of skew-symmetrizable matrices can be
extended to a natural equivalence of properly chosen \QC companions.
The results of these two sections are put together in
Section~\ref{sec:proof}, where the proof of Theorem \ref{thm:MainResult}
is completed.

The concluding Section~\ref{sec:cog} is purely graph-theoretic.
We call a graph~$\Gam$ \emph{cyclically orientable}
if it admits an orientation in which any chordless cycle is cyclically
oriented.
(For example, the full graph on four vertices is \emph{not} \CO.)
The main result of Section~\ref{sec:cog} is Theorem~\ref{thm:CO-MainResult}
which gives several properties of graphs that are equivalent to \COP.
As a consequence, we obtain a graph-theoretic statement
(Corollary~\ref{cor:product-signs}) which implies
Proposition~\ref{pr:cyclic-B-to-A}.
The paper concludes with Remark~\ref{rem:checking-4} discussing a
possible strategy for checking condition~(4).

\begin{remark}
\label{rmk:MainResult2}
The preliminary version of Theorem~\ref{thm:MainResult} due to the first
two authors was stated in terms of \emph{integral quadratic forms} rather
than \qCms. (Recall that the quadratic form represented by a symmetric matrix
$C$ with positive diagonal entries is called \emph{integral} if
$2C_{ij}/C_{ii}\in\ZZ$ for all $i,j$; in our present terminology,
this means that $C = DA$ is the symmetrized version of a \qCm $A$.)
The classification of positive definite integral quadratic forms
(up to natural equivalence) is well-known, and the answer is once again
given by the Cartan-Killing classification, cf.
Proposition~\ref{pr:positive-qCm-classification}.
Our results imply that this classification is in agreement with the one
given by Theorem~\ref{thm:FZ}: if~$A$ is a positive \QC companion of~$B$,
then the Cartan-Killing type of the cluster algebra associated
to~$B$ is the same as the Cartan-Killing type of the quadratic
form represented by the symmetrized matrix~$DA$.
\end{remark}

\section{Properties of positive \qCms}
\label{sec:pos_properties}
This section is a ``symmetrizable analogue'' of~\cite[Section~9]{FZ2}.
We start with a simple observation.

\begin{lemma}
\label{lem:2-finite}
Let $A$ be a positive \qCm. Then
\begin{itemize}
\item[(a)]
$0\leq A_{ij}A_{ji} \leq 3$ for any $i\neq j$.
\item[(b)]
$A_{ik}A_{kj}A_{ji}\geq 0$ for any pairwise different $i,j,k$.
\end{itemize}
\end{lemma}

\proof
(a) is immediate from the positivity of the principal minor
of~$A$ on the rows and columns~$i$ and~$j$.

(b) Let $A_{ik} A_{kj} A_{ji} \neq 0$.
Since $A$ is symmetrizable, we have
$A_{ki} A_{jk} A_{ij}  = A_{ik} A_{kj} A_{ji}$.
The positivity condition for the principal $3 \times 3$ minor of $A$
on the rows and columns $i, j, k$ can now be rewritten as
\begin{equation}
\label{eq:det}
A_{ik} A_{kj} A_{ji}  > A_{ij} A_{ji} + A_{ik} A_{ki} + A_{jk}
A_{kj} - 4 \geq -1,
\end{equation}
implying our claim.
\endproof

We now introduce the \emph{diagram} of a \qCm, which is a symmetrizable
analogue of~\cite[Definition~7.3]{FZ2}.

\begin{definition}
\label{def:A-diagram}
The \emph{diagram} $\Gam(A)$ of a
$n \times n$ \qCm~$A$ is a (undirected) graph with vertices
$\{1,2,\ldots, n\}$ and edges $\{i,j\}$ for each $i\neq j$ with
$A_{ij}\neq 0$, where every edge $\{i,j\}$ is assigned the \emph{weight}
$A_{ij} A_{ji}$ and the \emph{sign}
$\varepsilon_{ij} = -\sgn(A_{ij}) = -\sgn(A_{ji})$.
\end{definition}

In drawing the diagrams, all unspecified weights will be assumed
to be equal to~$1$.
With some abuse of notation, we denote by the same symbol
$\Gam(A)$ the underlying edge-weighted graph of the diagram, obtained by
forgetting signs of the edges, and also the underlying graph obtained by
forgetting both signs and edge weights.
Note that an edge-weighted graph~$\Gam$ which is of the form~$\Gam(A)$
must satisfy the following condition (see~\cite[Exercise~2.1]{kac}) :
\begin{equation}
\label{eq:perfect-square}
\text{The product of edge weights along every cycle
of~$\Gam$ is a perfect square.}
\end{equation}

\begin{definition}
\label{def:positive-graph}
An edge-weighted graph will be called \emph{positive}
if some sign assignment to the edges makes it into the diagram
$\Gam(A)$ of some positive \qCm~$A$.
\end{definition}

The following proposition is an analogue of~\cite[Proposition~9.3]{FZ2}.

\begin{proposition}
\label{prp:Dyn}
Positive edge-weighted trees are precisely Dynkin diagrams.
Each of them becomes the diagram of a positive \qCm under
an arbitrary assignment of signs.
\end{proposition}

\proof
Suppose $\Gam(A)$ is a tree.
We can assume without loss of generality that the signs of all
edges are equal to~$1$, i.e., $A$ is a generalized Cartan matrix
(this can be achieved by replacing $A$
if necessary by a positive \qCm of the form $A'= E A E$,
where $E$ is a diagonal matrix with entries $E_{ii}=\pm 1$).
Our statement now follows from the Cartan-Killing classification
of positive generalized Cartan matrices, see, e.g.,
\cite[Theorem~4.8]{kac} or Theorems~1 and~4 in~\cite[VI,4]{B}.
\endproof

In analogy with~\cite[Definition~9.1]{FZ2}, by a
\emph{subdiagram} of $\Gam(A)$ we mean a diagram of the form
$\Gam(A')$, where~$A'$ is a principal submatrix of~$A$.
Thus, $\Gam(A')$ is obtained from $\Gam(A)$ by taking an
induced subgraph on a subset of vertices and
keeping all the edge weights and signs the same as in $\Gam(A)$.
Since the positivity of~$A$ implies that of~$A'$, we obtain the
following corollary.

\begin{corollary}
\label{lem:extendedDynkin}
None of the following edge-weighted graphs can appear as subdiagrams of
the diagram of a positive \qCm:
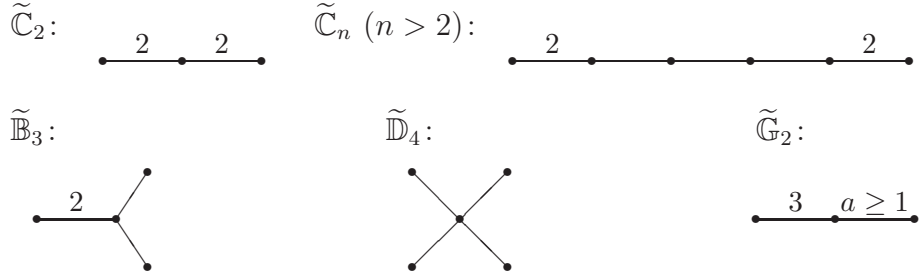
\begin{center}
  \begin{picture}(352,105)
    \put(0,80){
      \put(-10,10){$\widetilde{\mathbb{C}}_{2}\!:$}
      \multiput(25,0)(30,0){3}{\circle*{3}}
      \put(25,0){\line(1,0){60}}
      \multiput(40,3)(30,0){2}{\HBCenter{\small $2$}}
      }
    \put(120,80){
      \put(-15,10){$\widetilde{\mathbb{C}}_{n}$ $(n>2)\!:$}
      \multiput(60,0)(30,0){6}{\circle*{3}}
      \put(60,0){\line(1,0){150}}
      \multiput(75,3)(120,0){2}{\HBCenter{\small $2$}}
      }
    \put(0,20){
      \put(-10,30){$\widetilde{\mathbb{B}}_{3}\!:$}
      \multiput(0,0)(30,0){2}{\circle*{3}}
      \multiput(42,-18)(0,36){2}{\circle*{3}}
      \put(0,0){\line(1,0){30}}
      \put(30,0){\line(2,-3){12}}
      \put(30,0){\line(2,3){12}}
      \put(15,3){\HBCenter{\small $2$}}
      }
    \put(142,20){
      \put(-10,30){$\widetilde{\mathbb{D}}_{4}\!:$}
      \multiput(0,-18)(18,18){3}{\circle*{3}}
      \multiput(0,18)(36,-36){2}{\circle*{3}}
      \put(0,-18){\line(1,1){36}}
      \put(0,18){\line(1,-1){36}}
      }
    \put(272,20){
      \put(0,30){$\widetilde{\mathbb{G}}_{2}\!:$}
      \multiput(0,0)(30,0){3}{\circle*{3}}
      \put(0,0){\line(1,0){60}}
      \put(15,3){\HBCenter{\small $3$}}
      \put(45,3){\HBCenter{\small $a \geq 1$}}
      }
  \end{picture}
\end{center}
\end{corollary}

Our next result is an analogue of~\cite[Proposition~9.7]{FZ2}.

\begin{proposition}
\label{lem:positive_cycles}
Positive edge-weighted cycles are precisely those of the following
three types:
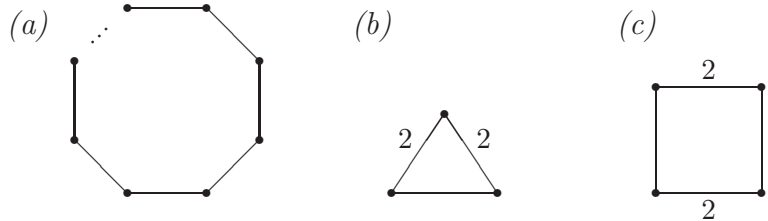
\begin{equation*}
\begin{picture}(260,90)
\put(0,10){
  \put(-10,60){\RBCenter{(a)}}
  \multiput(0,20)(0,30){2}{\circle*{3}}
  \multiput(70,20)(0,30){2}{\circle*{3}}
  \multiput(20,0)(30,0){2}{\circle*{3}}
  \multiput(20,70)(30,0){2}{\circle*{3}}
  \put(20,70){\line(1,0){30}}
  \put(50,70){\line(1,-1){20}}
  \put(70,50){\line(0,-1){30}}
  \put(70,20){\line(-1,-1){20}}
  \put(50,0){\line(-1,0){30}}
  \put(20,0){\line(-1,1){20}}
  \put(0,20){\line(0,1){30}}
  \multiput(7,57)(2.5,2.5){3}{\circle*{1}}
}
\put(120,10){
  \put(0,60){\RBCenter{(b)}}
  \multiput(0,0)(40,0){2}{\circle*{3}}
  \put(20,30){\circle*{3}}
  \put(0,0){\line(1,0){40}}
  \put(0,0){\line(2,3){20}}
  \put(40,0){\line(-2,3){20}}
  \put(8,17){\RBCenter{\small $2$}}
  \put(32,17){\small $2$}
}
\put(220,10){
  \put(0,60){\RBCenter{(c)}}
  \multiput(0,0)(40,0){2}{\circle*{3}}
  \multiput(0,40)(40,0){2}{\circle*{3}}
  \multiput(0,0)(40,0){2}{\line(0,1){40}}
  \multiput(0,0)(0,40){2}{\line(1,0){40}}
  \put(20,43){\HBCenter{\small $2$}}
  \put(20,-3){\HTCenter{\small $2$}}
}
\end{picture}
\end{equation*}
Furthermore, a sign assignment makes each of these cycles the
diagram of a positive \qCm if and only if the product of signs along all
the edges is equal to~$-1$.
\end{proposition}

\proof
First suppose that~$A$ is a positive \qCm whose diagram is a
cycle not of type~(a).
Then at least one edge of $\Gam(A)$ has weight $a>1$.
By Lemma~\ref{lem:2-finite}(a), the weight~$a$ is either~$2$ or $3$.
It follows from~\eqref{eq:perfect-square} that at least two edges of $\Gam(A)$ have weight~$a$.
If $a=3$ then the cycle must be a triangle since otherwise it
would have a subdiagram $\widetilde{\mathbb{G}}_2$,
a contradiction to Corollary~\ref{lem:extendedDynkin}.
We are left with the case of a $3 \times 3$ matrix~$A$
such that $A_{12}A_{21}=A_{23}A_{32} = 3$, and
$A_{13}A_{31}=1$.
In view of~\eqref{eq:det}, such a matrix~$A$ cannot be positive.

So, we can assume that all the edge weights are equal to~$1$ or~$2$,
with at least two edges of weight~$2$.
Then~$\Gam(A)$ must be of one of the types (b) or (c),
since otherwise it would have a subdiagram $\widetilde{\mathbb{C}}_{n}$
for some $n  \geq 2$, again in contradiction to
Corollary~\ref{lem:extendedDynkin}.

To finish the proof, it remains to show that, in each of the cases
(a), (b) and (c), the positivity of a \qCm~$A$ is equivalent to
the condition that the product of signs along all the edges is equal to~$-1$.
Since in each case, all the proper subdiagrams of~$\Gam(A)$ are
Dynkin diagrams, we only need to show that the above sign
condition is equivalent to $\det(A) > 0$.
Using simultaneous permutations and changes of signs of rows and columns,
we can assume that~$A_{ij} < 0$ for $|j-i| = 1$; and we need to show that
$\det(A) > 0$ precisely when $A_{1n}$ and $A_{n1}$ are positive.

Dealing with case (a) first, let $A_n(x)$ be the $n \times n$ \qCm
with non-zero off-diagonal entries $A_{ij} = -1$ for $|j-i| = 1$,
and $A_{1n}=A_{n1}=x$.
By a standard calculation,
$$\det(A_n(x)) = n+1 + 2x - x^2(n-1);$$
therefore, $\det(A_n(1)) = 4$, and $\det(A_n(-1)) = 0$, as required.
The cases (b) and (c) are similar (and simpler).
\endproof

\begin{remark}
\label{rem:sign-condition}
The sign condition in Proposition~\ref{lem:positive_cycles}
implies Proposition~\ref{pr:cyclic-A}.
\end{remark}

\begin{definition}
\label{def:qCm-equivalence}
We say that \qCms~$A$ and~$A'$ are \emph{equivalent} and write
$A' \sim A$ if $A$ and $A'$ have the same symmetrizer~$D$
(that is, $D$ is a diagonal matrix with positive diagonal entries such that
$C = DA$ and $C' = DA'$ are symmetric), and the symmetrized matrices satisfy
$C' = E^T C E$ for some integer matrix~$E$ with
determinant~$\pm 1$.
\end{definition}

Note that this definition does not depend on the choice of a symmetrizer~$D$.

We conclude this section by showing that the equivalence classes
of positive \qCms are classified by Cartan-Killing types.
Although this result is well-known to experts, we were unable to
find an adequate reference, so will outline the proofs.

Let $A$ be a $n \times n$ \qCm.
For each $i = 1, \dots, n$, define an automorphism $s_i$ of the lattice
$\ZZ^n$ by setting $s_i (e_j) = e_j - A_{ij} e_i$,
where $\{e_1, \dots, e_n\}$ is the standard basis in~$\ZZ^n$.
Let $W(A) \subset GL_n(\ZZ)$ be the group generated by
$s_1, \dots, s_n$.

\begin{proposition}
\label{pr:positive-qCm-classification}
The following conditions on a \qCm~$A$ are equivalent:
\begin{enumerate}
\item [\rm (1)]
$A$ is positive.

\item [\rm (2)]
The group $W(A)$ is finite.

\item [\rm (3)]
There exist a (reduced) root system $\Phi$
and a linearly independent subset $\{\beta_1, \dots, \beta_n\} \subset \Phi$
such that $A_{ij} = \langle \beta_i^\vee, \beta_j \rangle$,
where $\beta^\vee$ is the coroot dual to a root $\beta$.

\item [\rm (4)]
$A$ is equivalent to a positive Cartan matrix~$A^0$.
\end{enumerate}
Under these conditions, if $\Phi_0 \subset \Phi$ is the smallest
root subsystem of~$\Phi$ that contains the set $\{\beta_1, \dots, \beta_n\}$
in {\rm (3)}, then the Cartan-Killing type of~$\Phi_0$
is the same as the Cartan-Killing type
of the matrix~$A^0$ in {\rm (4)}, and it characterizes~$A$ up to
equivalence.
Furthermore,~$W(A)$ is naturally identified with the
Weyl group of~$\Phi_0$.
\end{proposition}

\proof
The implication $(4) \Longrightarrow (1)$ is trivial.
To prove $(1) \Longrightarrow (2)$, suppose that~$A$ is positive,
and let $(\alpha |\ \beta)$ be the positive definite (symmetric)
scalar product in $\ZZ^n$ given by $(e_i |\ e_j) = C_{ij}$,
where $C = (C_{ij}) = DA$ is the positive definite symmetric
matrix corresponding to~$A$.
The standard  check shows that, with respect to this scalar product,
$s_i$ is the orthogonal reflection in the
orthogonal complement to~$e_i$.
Thus $W$ is a discrete subgroup of the compact orthogonal group
$O_n(\RR)$, and so is finite, as claimed.

The implication $(2) \Longrightarrow (3)$ follows from a well-known
classification of finite crystallographic reflection groups
given in \cite{B}: namely, every such group is the Weyl
group of a reduced root system~$\Phi$.
To be more specific, take
$\Phi = W(A)\,\{e_1, \dots, e_n\} \subset \ZZ^n$.
If $W(A)$ is finite then so is $\Phi$; the fact that $\Phi$
satisfies the rest of the axioms of the root systems in \cite{B}
is checked easily.
Then {\rm (3)} holds with $\beta_i = e_i$.

It remains to prove $(3) \Longrightarrow (4)$.
Without loss of generality, we can assume that $\Phi$ is
the smallest root system containing $\{\beta_1, \dots, \beta_n\}$,
that is, $\Phi = W \{\beta_1, \dots, \beta_n\}$, where $W = W(A)$
is the group generated by the reflections corresponding to the roots
$\beta_1, \dots, \beta_n$.
By this assumption, $\Phi$ has rank~$n$, and $\{\beta_1, \dots, \beta_n\}$
is a $\ZZ$-basis of the root lattice of~$\Phi$.

Let $(\alpha |\ \beta)$ be a $W$-invariant positive definite scalar
product in the root lattice.
By the definition, $A_{ij} = 2(\beta_i |\ \beta_j)/(\beta_i |\ \beta_i)$.
Thus, the symmetric matrix $C = DA$ corresponding to~$A$ has entries
$C_{ij} = (\beta_i |\ \beta_j)$, with the diagonal entries of~$D$
given by $D_{ii} = (\beta_i |\ \beta_i)/2$.

Let $\{\alpha_1, \dots, \alpha_n\}$ be a system of simple
roots in~$\Phi$, and $A^0$ the corresponding (positive) Cartan
matrix, so that
$A^0_{ij} = \langle \alpha_i^\vee, \alpha_j \rangle$.
Thus, the symmetric matrix $C^0 = D^0 A^0$ corresponding to~$A^0$ has entries
$C^0_{ij} = (\alpha_i |\ \alpha_j)$, with the diagonal entries of~$D^0$
given by $D^0_{ii} = (\alpha_i |\ \alpha_i)/2$.
It follows that $C = E^T C^0 E$,
where $E$ is the transition matrix from the basis
$\{\alpha_1,\ldots\alpha_n\}$ to the basis $\{\beta_1,\ldots,\beta_n\}$
(since both families are $\ZZ$-bases of the root lattice, the
matrix~$E$ is invertible over the integers).
To finish the proof it remains to show that the simple roots
$\alpha_1, \dots, \alpha_n$ can be chosen in such a way that $D^0 = D$.
In other words, we need to show the following:
\begin{equation}
\label{eq:D=D0}
\text{One can reorder the roots $\beta_1, \dots, \beta_n$
so that $(\beta_i |\ \beta_i) = (\alpha_i |\ \alpha_i)$ for all~$i$.}
\end{equation}

We can assume without loss of generality that the root
system~$\Phi$ is irreducible.
If~$\Phi$ is simply-laced, i.e., all roots are of the same length,
\eqref{eq:D=D0} becomes tautological.
So it is enough to check \eqref{eq:D=D0} for~$\Phi$ of one of the types
$\mathbb{B}_n, \mathbb{C}_n, \mathbb{F}_4$ and $\mathbb{G}_2$.
In each case, there are two different root lengths (``long" and ``short").
Our assumption that $\Phi$ is the smallest root system
containing $\{\beta_1, \dots, \beta_n\}$ implies that among the
roots~$\beta_i$ there is a long one and a short one.
This finishes the story for~$\mathbb{G}_2$.

Now let~$\Phi$ be of type~$\mathbb{B}_n$, with the standard choice of
simple roots
$\{\alpha_1, \dots, \alpha_n\} = \{e_1 - e_2, \dots, e_{n-1} - e_n, e_n\}$,
where $e_1, \dots, e_n$ is the standard basis in the root lattice $\ZZ^n$.
The long (resp. short) roots form a subsystem of type $\mathbb{D}_n$
(resp.~$\mathbb{A}_1^n$) consisting of the roots
$\pm e_i \pm e_j \,\, (1 \leq i < j \leq n)$ (resp. $\pm e_i \,\,
(1 \leq i \leq n)$).
To prove \eqref{eq:D=D0}, we observe the following: every subset
of~$\Phi$ containing at most $n-2$ long roots, is contained in
a root subsystem of type $\mathbb{B}_k \times \mathbb{B}_{n-k}$ for some $k = 1, \dots, n-1$
(to see this, represent each long root $\pm e_i \pm e_j$ by an edge
joining~$i$ and~$j$, and use the obvious fact that every graph with~$n$ vertices and
at most $n-2$ edges is disconnected).
Therefore, the set $\{\beta_1, \dots, \beta_n\}$ must consist of
$n-1$ long roots and a short one, proving  \eqref{eq:D=D0}.
The type~$\mathbb{C}_n$ is treated in the same way (with long and short
roots interchanged).

Finally, let~$\Phi$ be of type~$\mathbb{F}_4$, with the simple roots
$\{\alpha_1, \dots, \alpha_4\} = \{e_2 - e_3, e_3 - e_4, e_4,
\frac{1}{2}(e_1 - e_2 - e_3 - e_4)\}$ (see \cite{B}).
In view of the symmetry between the long and short roots,
to prove \eqref{eq:D=D0} it suffices to show that
$\{\beta_1, \dots, \beta_4\}$ cannot consist
of one long root and three short ones.
Suppose on the contrary that, say $\beta_1 = \pm e_1 \pm e_2$,
and each of the roots $\beta_2, \beta_3$ and $\beta_4$ is of the
form $\pm e_i$ or $\frac{1}{2}(\pm e_1 \pm e_2 \pm e_3 \pm e_4)$.
An easy inspection shows that $W \{\beta_1, \dots, \beta_4\}$ is then
contained in the root subsystem of~$\Phi$ consisting of all short
roots and the eight long roots
$\{\pm e_1 \pm e_2, \pm e_3 \pm e_4\}$, again in contradiction with
our assumption that $\Phi$ is the smallest root system
containing $\{\beta_1, \dots, \beta_4\}$.

This concludes the proof of \eqref{eq:D=D0} and hence of the
equivalence of conditions $(1) - (4)$.
It is well known that positive Cartan matrices of different types
are not equivalent to each other.
The rest of the statements in Proposition~\ref{pr:positive-qCm-classification}
have been already established in the course of the above argument.
\endproof

\begin{remark}
\label{rem:Ringel-unit-forms}
In the case where a \qCm is assumed to be \emph{symmetric},
a proof of Proposition~\ref{pr:positive-qCm-classification} was essentially given
in~\cite[1.2]{Ri} (using the language of unit quadratic forms).
\end{remark}

\section{Mutations and \QC companions}
\label{sec:mutations}

In this section, we show that, under some conditions, the
mutation-equivalence of skew-symmetrizable matrices can be
extended to an equivalence (in the sense of Definition~\ref{def:qCm-equivalence})
of properly chosen \QC companions.

\begin{definition}
\label{def:k-compatible}
For a matrix index~$k$, we say that a
\QC companion~$A$ of a skew-symmetrizable matrix~$B$ is
\emph{$k$-compatible} with~$B$
if the signs of its entries satisfy
\begin{equation}
\label{eq:k-compatible}
\text{if $B_{ik} > 0$ and $B_{kj} > 0$ for some $i, j$
then $\sgn(A_{ik} A_{kj} A_{ji}) = \sgn(B_{ji})$.}
\end{equation}
\end{definition}

\begin{proposition}
\label{pr:mutation-S-counterpart}
Let~$A$ be a $k$-compatible \QC companion of a
skew-sym\-met\-rizable matrix~$B$.
Then there exists a $k$-compatible \QC companion~$A'$ of
$B' = \mu_k(B)$ such that~$A' \sim A$.
Explicitly, the signs of off-diagonal
matrix entries of~$A'$ can be chosen as follows:
\begin{equation}
\label{eq:sgn-A'}
\sgn(A'_{ij}) =\sgn(A'_{ji})=
\begin{cases}
\sgn(B_{ik}) \sgn(A_{ik})& \text{if $i \neq k, \, j = k$;} \\
\sgn(B_{ij} B'_{ij}) \sgn(A_{ij})& \text{if $i \neq k, \, j \neq k$.}
\end{cases}
\end{equation}
\end{proposition}

\proof
We fix~$k$ and introduce the following three $n \times n$
matrices:
\begin{itemize}
\item $J$ is the diagonal matrix with $J_{kk}=-1$ and $J_{ii}=1$ for
$i \neq k$.
\item $E$ is the matrix with all the entries outside the $k$-th
column equal to~$0$, and the $k$-th column entries given by
\begin{equation*}
E_{ik} =
\begin{cases}
A_{ik} & \text{if $B_{ik} > 0$;} \\
0 & \text{otherwise.}
\end{cases}
\end{equation*}
\item $F$ is the matrix with all the entries outside the $k$-th
row equal to~$0$, and the $k$-th row entries given by
\begin{equation*}
F_{kj} =
\begin{cases}
A_{kj} & \text{if $B_{kj} < 0$;} \\
0 & \text{otherwise.}
\end{cases}
\end{equation*}
\end{itemize}

We now set
\begin{equation}
\label{eq:A'-product}
A' = (J - E) A (J - F) \, ,
\end{equation}
and claim that~$A'$ satisfies all the required properties.
First, a direct calculation shows that the entries of
$A' = JAJ - EAJ - JAF + EAF$ are given by
\begin{equation}
\label{eq:A'-entries}
A'_{ij} =
\begin{cases}
2 & \text{if $i = j = k$;}\\
\sgn(B_{ik}) A_{ik}& \text{if $i \neq k = j$;} \\
-\sgn(B_{kj}) A_{kj}& \text{if $i = k \neq j$;} \\
A_{ij} - \sgn(A_{ik} A_{kj}) [B_{ik} B_{kj}]_+&
\text{if $i \neq k, \, j \neq k$.}
\end{cases}
\end{equation}
Comparing \eqref{eq:A'-product} with \eqref{eq:matrix-mutation}
and using \eqref{eq:k-compatible}, it is easy to see that~$A'$ is
a \QC companion of~$B'$, and that it satisfies \eqref{eq:sgn-A'}.
Furthermore, the claim that~$A'$ is $k$-compatible with~$B'$ is a
direct consequence of \eqref{eq:sgn-A'}.
Finally, to show that $A' \sim A$, we note that
$$D(J - E) = (D(J - F))^T = (J - F)^T D.$$
Therefore,
$$C' = DA' = D (J - E) A (J - F) = (J - F)^T D A (J - F) =
(J - F)^T C (J - F) \, ,$$
and we are done.
\endproof

\begin{corollary}
\label{cor:compatible-positive}
Suppose a skew-symmetrizable matrix~$B$ satisfies condition~(4) in
Theorem~\ref{thm:MainResult}, and let~$A$ be a positive \QC companion of~$B$.
Then~$\mu_k(B)$ has a positive \QC companion for any matrix index~$k$.
\end{corollary}

\proof
In view of Lemma~\ref{lem:2-finite}~(b),~$A$ is $k$-compatible with~$B$.
By Proposition~\ref{pr:mutation-S-counterpart}, the matrix~$\mu_k(B)$
has a \QC companion which is equivalent to~$A$ and so is positive.
\endproof

In the rest of this section, we use Proposition~\ref{pr:mutation-S-counterpart}
to obtain more information on positive edge-weighted graphs, see
Definition~\ref{def:positive-graph}.
First a piece of notation.
Let~$\Gam$ be an edge-weighted graph,~$i \neq j$ two vertices
of~$\Gam$, and~$t$ a positive integer.
We denote by $\Gam[i,j;t]$ the edge-weighted graph obtained
from~$\Gam$ by adjoining~$t$ new vertices~$i_1, \dots, i_t$ and
$t+1$ new edges $\{i,i_1\}, \{i_1, i_2\}, \dots, \{i_{t-1}, i_t\}, \{i_t, j\}$,
all of weight one.

\begin{lemma}[Chain Contraction]
\label{lem:chain-contraction}
If the graph $\Gam[i,j;t]$ is positive then so is
$\Gam[i,j;1]$.
\end{lemma}

As a preparation for the proof, we recall several facts from \cite{FZ2}.
First recall from~\cite[Definition~7.3]{FZ2} that every
skew-symmetrizable matrix~$B$ has the \emph{diagram} $\Gam(B)$
which is a directed edge-weighted graph.
As an edge-weighted graph, it is completely analogous to its
symmetrizable counterpart in Definition~\ref{def:A-diagram}:
the vertices correspond to matrix indices, with an edge $\{i,j\}$ for each
$i\neq j$ with $B_{ij}\neq 0$, supplied with the \emph{weight}
$|B_{ij} B_{ji}|$.
Note that the edge weights satisfy~\eqref{eq:perfect-square}.
The only difference between the two kinds of diagrams is in the
way of encoding the signs of matrix entries: while for a \qCm~$A$,
an edge $\{i,j\}$ in $\Gam(A)$ is assigned the
sign $\varepsilon_{ij} = -\sgn(A_{ij}) = -\sgn(A_{ji})$, in
$\Gam(B)$ we have an arrow $i \to j$ whenever $B_{ij} > 0$
(and so $B_{ji} < 0$).
As before, in drawing the diagram $\Gam(B)$,
all unspecified weights will be equal to~$1$.

The next proposition reproduces \cite[Proposition~8.1]{FZ2}.

\begin{figure}[ht]
\setlength{\unitlength}{1.5pt}
\begin{picture}(30,17)(-5,0)
\put(0,0){\line(1,0){20}}
\put(0,0){\line(2,3){10}}
\put(0,0){\vector(2,3){6}}
\put(10,15){\line(2,-3){10}}
\put(10,15){\vector(2,-3){6}}
\put(0,0){\circle*{2}}
\put(20,0){\circle*{2}}
\put(10,15){\circle*{2}}
\put(2,10){\makebox(0,0){$a$}}
\put(18,10){\makebox(0,0){$b$}}
\put(10,-4){\makebox(0,0){$c$}}
\put(10,19){\makebox(0,0){$k$}}
\end{picture}
$
\begin{array}{c}
\stackrel{\textstyle\mu_k}{\longleftrightarrow}
\\[.3in]
\end{array}
$
\setlength{\unitlength}{1.5pt}
\begin{picture}(30,17)(-5,0)
\put(0,0){\line(1,0){20}}
\put(0,0){\line(2,3){10}}
\put(10,15){\vector(-2,-3){6}}
\put(10,15){\line(2,-3){10}}
\put(20,0){\vector(-2,3){6}}
\put(0,0){\circle*{2}}
\put(20,0){\circle*{2}}
\put(10,15){\circle*{2}}
\put(2,10){\makebox(0,0){$a$}}
\put(18,10){\makebox(0,0){$b$}}
\put(10,-4){\makebox(0,0){$c'$}}
\put(10,19){\makebox(0,0){$k$}}
\end{picture}

\vspace{-.2in}
\caption{Diagram mutation}
\label{fig:diagram-mutation-general}
\end{figure}
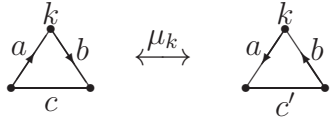

\begin{proposition}
\label{pr:diagram-mutation}
For a skew-symmetrizable matrix $B$,
the diagram $\Gam'\!\!=\!\Gam(\mu_k(B))$ is uniquely determined by
the diagram $\Gam=\Gam(B)$ and a matrix index~$k$.
Specifically, $\Gam'$ is obtained from $\Gam$ as follows:
\begin{itemize}
\item The orientations of all edges incident to~$k$ are reversed,
their weights intact.
\item
For any vertices $i$ and $j$ which are connected in
$\Gam$ via a two-edge oriented path going through~$k$ (refer to
Figure~\ref{fig:diagram-mutation-general} for the rest of notation),
the direction of the edge $(i,j)$ in $\Gam'$ and its weight $c'$
are uniquely determined by the rule
\begin{equation}
\label{eq:weight-relation-general}
\pm\sqrt {c} \pm\sqrt {c'} = \sqrt {ab} \,,
\end{equation}
where the sign before $\sqrt {c}$
(resp., before $\sqrt {c'}$)
is ``$+$'' if $i,j,k$ form an oriented cycle
in~$\Gam$ (resp., in~$\Gam'$), and is ``$-$'' otherwise.
Here either $c$ or $c'$ can be equal to~$0$.

\item
The rest of the edges and their weights in $\Gam$
remain unchanged.
\end{itemize}
\end{proposition}

\noindent{\bf Proof of Lemma~\ref{lem:chain-contraction}.}
It suffices to show the following:
\begin{equation}
\label{eq:t-t-1}
\text{If $t \geq 2$ then the positivity of $\Gam[i,j;t]$ implies
that of $\Gam[i,j;t-1]$.}
\end{equation}
Choose any skew-symmetrizable matrix~$B$ such that
$\Gam(B) = \Gam[i,j;t]$ (as edge-weighted graphs), and the edges
$\{i,i_1\}, \{i_1, i_2\}, \dots, \{i_{t-1}, i_t\}, \{i_t, j\}$ are
oriented as $i \to i_1 \to \cdots \to i_t \to j$.
Since $\Gam[i,j;t]$ is assumed to be positive,~$B$ has a positive
\QC companion~$A$.
Clearly, our choice of orientations implies that~$A$ is
$i_t$-compatible with~$B$.
By Proposition~\ref{pr:mutation-S-counterpart}, the diagram
$\Gam' = \Gam(\mu_{i_t}(B))$ is positive (as an edge-weighted graph).
By \eqref{eq:weight-relation-general}, $\Gam'$ is obtained from
$\Gam[i,j;t]$ by adjoining the edge $\{i_{t-1},j\}$.
Thus, $\Gam'$ has $\Gam[i,j;t-1]$ as an induced subgraph, and so
$\Gam[i,j;t-1]$ is positive.
This concludes the proof of \eqref{eq:t-t-1} and
Lemma~\ref{lem:chain-contraction}.
\endproof

We conclude this section with one more application of Lemma~\ref{lem:chain-contraction},
which allows us to construct a family of non-positive edge-weighted graphs.

\begin{lemma}
\label{lem:wheels-non-positive}
Let~$\Gam$ be an edge-weighted graph obtained from a cycle~$Z$
with all edge weights~$1$ by adjoining a vertex~$k$ joined with
an even number of vertices of~$Z$ with edges of weights not
exceeding~$2$.
Then~$\Gam$ is non-positive, with the only exception,
where~$k$ is joined with just two adjacent vertices of~$Z$,
both edges having weight~$1$.
\end{lemma}

\proof
First assume that at least one edge through~$k$ has weight~$2$.
In view of~\eqref{eq:perfect-square}, then every edge through~$k$
must have weight~$2$.
If two of these edges connect~$k$ with non-adjacent vertices
of~$Z$, then they form a subdiagram~$\widetilde{\mathbb{C}}_{2}$,
so~$\Gam$ cannot be positive.
This leaves us with the case where~$k$ is joined with just two adjacent
vertices of~$Z$, both edges having weight~$2$.
Then the graph $\Gam$ is of the form $\Delta[i,j;t]$, where
$\Delta$ is a triangle with edge weights~$2, 2$ and $1$, the edge
$\{i,j\}$ being of weight~$1$.
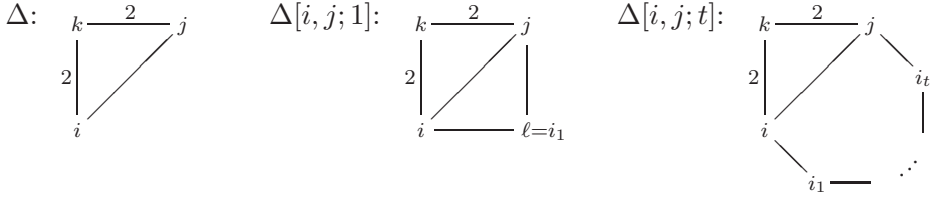
\begin{center}
\begin{picture}(310,82)(0,-27)
\put(0,0){
  \put(-15,40){\RBCenter{\small $\Delta$:}}
  \put(-2,0){\HVCenter{ $\scriptstyle i$}}
  \put(40,40){\HVCenter{$\scriptstyle j$}}
  \put(0,40){\HVCenter{$\scriptstyle k$}}
  \put(4,40){\line(1,0){31}}
  \put(0,6){\line(0,1){28}}
  \put(4,4){\line(1,1){32}}
  \put(-2,20){\RVCenter{$\scriptstyle 2$}}
  \put(20,42){\HBCenter{$\scriptstyle 2$}}
}
\put(130,0){
  \put(-15,40){\RBCenter{\small $\Delta[i,j;1]$:}}
  \put(0,0){\HVCenter{$\scriptstyle i$}}
  \put(40,40){\HVCenter{$\scriptstyle j$}}
  \put(38,0){\LVCenter{$\scriptstyle \ell=i_1$}}
  \put(0,40){\HVCenter{$\scriptstyle k$}}
  \put(5,0){\line(1,0){30}}
  \put(4,40){\line(1,0){31}}
  \put(0,6){\line(0,1){28}}
  \put(40,6){\line(0,1){26}}
  \put(4,4){\line(1,1){32}}
  \put(-2,20){\RVCenter{$\scriptstyle 2$}}
  \put(20,42){\HBCenter{$\scriptstyle 2$}}
}
\put(260,0){
  \put(-15,40){\RBCenter{\small $\Delta[i,j;t]$:}}
  \put(0,0){\HVCenter{$\scriptstyle i$}}
  \put(40,40){\HVCenter{$\scriptstyle j$}}
  \put(20,-20){\HVCenter{$\scriptstyle i_1$}}
  \put(60,20){\HVCenter{$\scriptstyle i_t$}}
  \put(0,40){\HVCenter{$\scriptstyle k$}}
  \put(4,-4){\line(1,-1){12}}
  \put(25,-20){\line(1,0){15}}
  \put(4,40){\line(1,0){31}}
  \put(0,6){\line(0,1){28}}
  \put(44,36){\line(1,-1){12}}
  \put(4,4){\line(1,1){32}}
  \put(-2,20){\RVCenter{$\scriptstyle 2$}}
  \put(20,42){\HBCenter{$\scriptstyle 2$}}
  \put(60,14){\line(0,-1){15}}
  \put(52,-17){\circle*{1}}
  \put(57,-12){\circle*{1}}
  \put(55,-15){\circle*{1}}
}
\end{picture}
\end{center}

By Lemma~\ref{lem:chain-contraction}, the non-positivity of~$\Gam$
follows from that of a $4$-vertex graph $\Delta[i,j;1]$ obtained
from~$\Delta$ by adjoining one more vertex~$\ell=i_1$ and two edges
$\{i,\ell\}$ and $\{j,\ell\}$, both of weight~$1$, as shown in the picture
above.
Arguing as in the proof of Lemma~\ref{lem:chain-contraction},
choose a skew-symmetrizable matrix~$B$ such that
$\Gam(B) = \Delta[i,j;1]$ (as edge-weighted graphs), and
both triangles $\{i,j,k\}$ and $\{i,j,\ell\}$ are cyclically
oriented.
Performing the mutation~$\mu_\ell$ destroys the edge~$\{i,j\}$,
and so transforms $\Delta[i,j;1]$ into a $4$-cycle with edge
weights (in cyclical order) $2, 2, 1, 1$.
By Proposition~\ref{lem:positive_cycles}, the latter cycle is
non-positive, hence by Corollary~\ref{cor:compatible-positive},
so are $\Delta[i,j;1]$ and $\Gam$.

It remains to consider the case where all the~$2p$ edges through~$k$
have weight~$1$.
We start with the case $p = 1$, that is, $k$ is joined with two
non-adjacent vertices of~$Z$.
Then~$\Gam$ has exactly three chordless cycles, with each edge
belonging to exactly two of these cycles.
This makes it impossible to attach the signs to edges to satisfy
the sign condition in Proposition~\ref{lem:positive_cycles},
so~$\Gam$ cannot be positive.

Next let~$p = 2$, that is, $k$ is joined with four vertices of~$Z$.
Applying Lemma~\ref{lem:chain-contraction} in the same way as
above, we can assume that~$Z$ has no other vertices, so~$\Gam$
is made up of four triangles.
Orienting all the edges so that these triangles become all
cyclically oriented, and performing the mutation~$\mu_k$
transforms~$\Gam$ into the graph~$\widetilde{\mathbb{D}}_{4}$.
Thus,~$\Gam$ is non-positive in this case as well.

If~$p=3$, arguing as above we can assume that~$\Gam$ is built from
a $6$-cycle~$Z$ with all its vertices joined with~$k$.
Orienting the edges of~$\Gamma$ so that these triangles become all
cyclically oriented, and performing mutations at two opposite
vertices of~$Z$ destroys four out of six edges through~$k$,
leaving us with the case~$p=1$ which we already dealt with.

Finally, if~$p \geq 4$ then $\Gam$ contains a
subdiagram~$\widetilde{\mathbb{D}}_{4}$, and so is again non-positive by
Corollary~\ref{lem:extendedDynkin}.
\endproof

\section{Proof of Theorem~\ref{thm:MainResult}}
\label{sec:proof}
The following result is the key ingredient of our proof of
Theorem~\ref{thm:MainResult}.

\begin{lemma}
\label{lem:mutation_closed}
Mutations of skew-symmetrizable matrices preserve
property {\rm (4)} in Theorem~\ref{thm:MainResult}.
\end{lemma}

\noindent{\bf Proof of Theorem \ref{thm:MainResult}.}
To deduce Theorem~\ref{thm:MainResult} from
Lemma~\ref{lem:mutation_closed}, we show that the conditions
$(2)$ and $(3)$ in Theorem~\ref{thm:FZ} satisfy the implications
$(2) \Longrightarrow (4) \Longrightarrow (3)$.
The implication $(4) \Longrightarrow (3)$ is immediate from
Lemma~\ref{lem:mutation_closed} and
Lemma~\ref{lem:2-finite}.
As for $(2) \Longrightarrow (4)$, it is again immediate from
Lemma~\ref{lem:mutation_closed} once we observe that
a matrix~$B$ in condition (2) of Theorem \ref{thm:FZ}
satisfies~(4) (by virtue of the Cartan-Killing classification,
the graph~$\Gamma(B)$ is a Dynkin diagram, hence a tree, so there
are no (chordless) cycles, and the corresponding condition in~(4)
is vacuous).\endproof

\smallskip

Turning to the proof of Lemma~\ref{lem:mutation_closed},
we will deduce it from a sharper statement
(Lemma~\ref{lem:quasi-finite-oriented} below).
To state it, we need some more terminology.

The following definition is motivated by
Propositions~\ref{prp:Dyn} and \ref{lem:positive_cycles}, and
Lemma~\ref{lem:wheels-non-positive}.

\begin{definition}
\label{def:quasi-finite}
An edge-weighted graph~$\Gam$
will be called \emph{quasi-finite}
if it satisfies the following conditions:
\begin{itemize}
\item
any induced subgraph of $\Gam$ which is a tree, is a Dynkin diagram.
\item
any induced subgraph of $\Gam$ which is a chordless cycle,
is of one of the types (a), (b), (c) in Proposition~\ref{lem:positive_cycles}.
\item
none of the non-positive edge-weight graphs in
Lemma~\ref{lem:wheels-non-positive} appears as an induced subgraph of~$\Gam$.
\end{itemize}
We say that a skew-symmetrizable matrix~$B$ and its diagram
$\Gam(B)$ are \QF if so is the underlying edge-weighted graph.
\end{definition}

\begin{remark}
\label{rmk:pos-are-QF}
Any skew-symmetrizable matrix $B$ which has a positive \QC companion is
\QF as follows from Propositions~\ref{prp:Dyn}
and~\ref{lem:positive_cycles}, and Lemma~\ref{lem:wheels-non-positive}.
\end{remark}

We will say that a
skew-symmetrizable matrix~$B$ (and its diagram $\Gam(B)$) is \COD if
every chordless cycle in $\Gam(B)$ is cyclically oriented.
We will deduce Lemma~\ref{lem:mutation_closed} from the following
statement.

\begin{lemma}
\label{lem:quasi-finite-oriented}
Suppose that two skew-symmetrizable matrices $B$ and $B'=\mu_k(B)$
are both \QF. Then $B$ is \COD if and only if so is $B'$.
\end{lemma}

Assuming Lemma 4.4, we can prove Lemma 4.1 as follows.

\noindent{\bf Proof of Lemma~\ref{lem:mutation_closed}.}
Suppose that a skew-symmetrizable matrix~$B$ satisfies~(4),
and let~$A$ be a positive \QC companion of~$B$.
We need to show that $B' = \mu_k(B)$  satisfies~(4) for any~$k$.
Note that~$B$ is \QF by Remark \ref{rmk:pos-are-QF} and \COD by the
definition of (4).
Note also that~$A$ is $k$-compatible with~$B$: to
check~\eqref{eq:k-compatible}, combine the fact that~$B$ is \COD
with Lemma~\ref{lem:2-finite}~(b).
Applying Proposition~\ref{pr:mutation-S-counterpart}, we see that~$B'
= \mu_k(B)$
has a positive \QC companion~$A'$ and therefore $B'$ is \QF, by Remark
\ref{rmk:pos-are-QF}.
By Lemma \ref{lem:quasi-finite-oriented}, $B'$ is \COD, and so
satisfies (4), as desired. \endproof

To finish the proof of Theorem~\ref{thm:MainResult}, it remains to
prove Lemma~\ref{lem:quasi-finite-oriented}.
Using Proposition~\ref{pr:diagram-mutation}, in our proof of
Lemma~\ref{lem:quasi-finite-oriented} we can forget about
skew-symmetrizable matrices and just work with
diagrams and their mutations given by the rules in
Proposition~\ref{pr:diagram-mutation}.
So in the following argument by a diagram we will just mean a
directed edge-weighted graph~$\Gam$ satisfying \eqref{eq:perfect-square}.

\medskip

\noindent{\bf Proof of Lemma~\ref{lem:quasi-finite-oriented}.}
Let~$\Gam$ be a diagram which is \QF but not \COD, and let
$\Gam'$ be obtained from $\Gam$ by mutation in some direction~$k$
(here we think of~$k$ as a vertex of~$\Gam$).
We need to show that~$\Gam'$ is either not \QF, or not \COD.
The proof will split into several cases.

\smallskip

\noindent{\bf Case 1.}~$k$ belongs to some chordless
cycle in~$\Gam$ which is not cyclically oriented.
Without loss of generality we can assume that this cycle is the
whole diagram~$\Gam$.
Let $\{i,k\}$ and $\{j,k\}$ be the two edges through~$k$ in~$\Gam$.

\noindent{\bf Case 1.1.}
The edges $\{i,k\}$ and $\{j,k\}$ are oriented both towards~$k$ or both
away from~$k$.
Then $\Gam'$ is obtained from~$\Gam$ by just reversing the orientations of
these edges, and so is not cyclically oriented.

\noindent{\bf Case 1.2.} The edges $\{i,k\}$ and $\{j,k\}$ are
oriented as $i \to k \to j$.

\noindent{\bf Case 1.2.1.} The cycle~$\Gam$ has at least four vertices.
Then in $\Gam'$ we have $j \to k \to i \to j$, and so, removing the vertex~$k$
from~$\Gam'$ leaves us with a cycle which is again not
cyclically oriented.

\noindent{\bf Case 1.2.2.} The cycle~$\Gam$ is a triangle with
vertices $i,j,k$ and weighted edges oriented as follows:
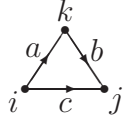
\begin{center}
\setlength{\unitlength}{1.5pt}
\begin{picture}(30,27)(-5,-5)
\put(0,0){\line(1,0){20}}
\put(0,0){\line(2,3){10}}
\put(0,0){\vector(2,3){6}}
\put(10,15){\line(2,-3){10}}
\put(10,15){\vector(2,-3){6}}
\put(7,0){\vector(1,0){6}}
\put(0,0){\circle*{2}}
\put(20,0){\circle*{2}}
\put(10,15){\circle*{2}}
\put(2,10){\makebox(0,0){$a$}}
\put(18,10){\makebox(0,0){$b$}}
\put(10,-4){\makebox(0,0){$c$}}
\put(-3,-3){\makebox(0,0){$i$}}
\put(23,-3){\makebox(0,0){$j$}}
\put(10,20){\makebox(0,0){$k$}}
\end{picture}
\end{center}

Applying \eqref{eq:weight-relation-general}, we see that the
weight~$c'$ of the edge $\{i,j\}$ in~$\Gam'$ satisfies
$$c' = (\sqrt{ab} + \sqrt{c})^2 = ab + 2 \sqrt{abc} + c \geq 4,$$
so in this case $\Gam'$ is not \QF.

\smallskip

\noindent{\bf Case 2.}~$k$ does not belong to any
non-oriented  chordless cycle in~$\Gam$.
Since~$\Gam$ is not \COD, it has at least one non-oriented
chordless cycle~$Z$.
Without loss of generality we can assume that $\Gam$ is just the
union of~$Z$ and~$\{k\}$.
In order for the mutation~$\mu_k$ to change~$Z$, there must be at least
two edges through~$k$ in~$\Gam$.
Furthermore, if we list these edges in the cyclical order
along~$Z$, then the incoming and outgoing edges must alternate, in
particular, there must be an even number of them; otherwise,~$k$
would belong to a non-oriented cycle.
Since~$\Gam$ is assumed to be \QF, the cycle~$Z$ is of one of the three
types (a), (b), or (c) in Proposition~\ref{lem:positive_cycles}.

\noindent{\bf Case 2.1.}~$Z$ is of type~(b), i.e., it is a
non-oriented triangle with edge weights $2$, $2$, and $1$.
Let~$k$ be joined with vertices~$i$ and~$j$ of~$Z$, with
orientations $i \to k \to j$.
Since the triangle $\{i,j,k\}$ in $\Gam$ must be oriented, the
edge $\{i,j\}$ is oriented as $j \to i$.
The mutation~$\mu_k$ reverses the orientations of the edges
$\{i,k\}$ and $\{j,k\}$, and changes the weight and orientation of the edge
$\{i,j\}$ in accordance with \eqref{eq:weight-relation-general}
(in particular, the latter edge can be absent in~$\Gam'$).
An easy inspection shows that the only chance for~$\Gam'$ to be
\COD, while~$\Gam$ is not \COD, is to have the edge $\{i,j\}$ to be of
opposite orientations in~$\Gam$ and $\Gam'$.
By~\eqref{eq:weight-relation-general}, this can only happen if
$\{i,j\}$ has weight~$1$ in~$\Gam$, while each of the remaining
four edges has weight~$2$.
But then~$\Gam$ has a subdiagram of type~$\widetilde{\mathbb{C}}_{2}$
(see Corollary~\ref{lem:extendedDynkin}), so cannot be \QF.

\noindent{\bf Case 2.2.}~$Z$ is of type~(c), i.e., it is a
non-oriented $4$-vertex cycle with edge weights (in cyclic order)
$1$, $2$, $1$, $2$.

\begin{figure}[ht]
\begin{picture}(332,110)
\put(0,5){
  \put(0,0){
    \put(36,88){\HBCenter{\small Case 2.2.1}}
    \put(72,0){\HVCenter{\small $i$}}
    \put(72,74){\HVCenter{\small $j$}}
    \put(36,36){\HVCenter{\small $k$}}
    \put(0,0){\line(0,1){72}}
    \put(72,6){\line(0,1){60}}
    \multiput(0,0)(0,72){2}{\line(1,0){66}}
    \multiput(0,0)(0,72){2}{\circle*{3}}
    \put(40,40){\line(1,1){28}}
    \put(40,32){\line(1,-1){28}}
    \multiput(50,16)(0,40){2}{\HVCenter{\small $2$}}
    \multiput(36,-5)(0,82){2}{\HVCenter{\small $2$}}
    }
  \put(130,0){
    \put(36,88){\HBCenter{\small Case 2.2.2}}
    \put(0,74){\HVCenter{\small $j$}}
    \put(72,0){\HVCenter{\small $i$}}
    \put(36,36){\HVCenter{\small $k$}}
    \multiput(0,0)(72,72){2}{\circle*{3}}
    \multiput(0,0)(72,6){2}{\line(0,1){66}}
    \multiput(0,0)(6,72){2}{\line(1,0){66}}
    \multiput(4,68)(36,-36){2}{\line(1,-1){28}}
    \put(32,40){\line(-1,1){28}}
    \put(16,50){\HVCenter{\small $2$}}
    \multiput(36,-5)(0,82){2}{\HVCenter{\small $2$}}
    }
  \put(260,0){
    \put(36,88){\HBCenter{\small Case 2.2.3}}
    \put(36,36){\HVCenter{\small $k$}}
    \multiput(0,0)(70,0){2}{\line(0,1){70}}
    \multiput(0,0)(0,70){2}{\line(1,0){70}}
    \multiput(0,70)(40,-40){2}{\line(1,-1){30}}
    \multiput(0,0)(40,40){2}{\line(1,1){30}}
    \multiput(49,16)(0,40){2}{\HVCenter{\small $2$}}
    \multiput(36,-5)(0,82){2}{\HVCenter{\small $2$}}
  }
}
\end{picture}
\caption{Different cases to be considered}
\label{fig:Cases-c}
\end{figure}
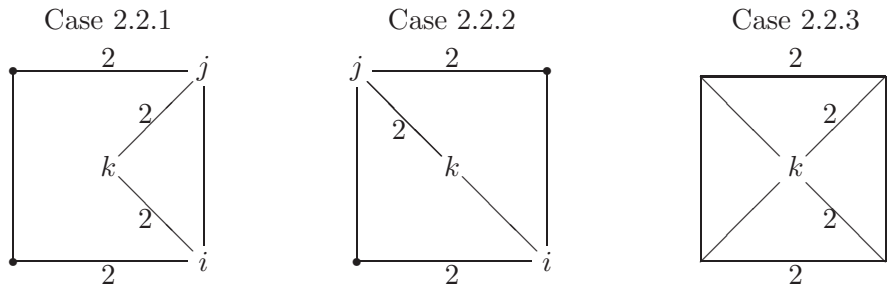

\noindent{\bf Case 2.2.1.}~$k$ is joined precisely with the two
adjacent vertices~$i$ and~$j$ of~$Z$.
The same argument as in Case~2.1 shows that, the only chance to
have~$\Gam'$ \COD and $\Gam$ not \COD appears if $\{i,j\}$ has
weight~$1$ in $\Gam$, while both $\{i,k\}$ and $\{j,k\}$ have weight~$2$.
But then~$\Gam$ has a subdiagram of
type~$\widetilde{\mathbb{C}}_{2}$, so cannot be \QF.

\noindent{\bf Case 2.2.2.}~$k$ is joined precisely with the two
opposite vertices~$i$ and~$j$ of~$Z$.
This case is even easier than the previous one: in view of
Proposition~\ref{lem:positive_cycles}, to satisfy
\eqref{eq:perfect-square}, the edges $\{i,k\}$ and $\{j,k\}$ must have weights~$1$ and~$2$,
which again forces~$\Gam$ to have a subdiagram of
type~$\widetilde{\mathbb{C}}_{2}$.

\noindent{\bf Case 2.2.3.}~$k$ is joined with all four vertices of~$Z$.
The only balance of edge weights that makes~$\Gam$ \QF is the
one shown in Figure \ref{fig:Cases-c}.
Furthermore, since we are still in Case~2, all four triangles
containing~$k$ must be cyclically oriented.
Performing the mutation~$\mu_k$, we obtain a diagram, which
contains a subdiagram $\widetilde{\mathbb{B}}_{3}$
and so is not \QF.

\noindent{\bf Case 2.3.}~$Z$ is of type~(a), i.e., has all
edge weights equal to~$1$.
According to Lemma~\ref{lem:wheels-non-positive}, if~$\Gam$ is \QF
then~$k$ is joined with just two adjacent vertices~$i$ and~$j$ of~$Z$,
both edges $\{i,k\}$ and $\{j,k\}$ having weight~$1$.
As before, the triangle $\{i,j,k\}$ must be cyclically oriented.
Applying the mutation~$\mu_k$ destroys the edge $\{i,j\}$, and so
transforms~$\Gam$ into a chordless cycle.
If this cycle is cyclically oriented, then so is the original cycle~$Z$,
contradicting our assumption.
This completes the proofs of Lemma~\ref{lem:quasi-finite-oriented}
and Theorem~\ref{thm:MainResult}.
\endproof

\section{Cyclically orientable graphs}
\label{sec:cog}

This section is purely graph-theoretic.
We call a graph~$\Gam$ \emph{cyclically orientable}
if it admits an orientation in which any chordless cycle is cyclically
oriented. For example, the full graph on four vertices is \emph{not} \CO.
We will give several properties of graphs that are equivalent to \COP.
This requires some notation.

For a given graph~$\Gam$, we denote by~$\Ver(\Gam)$ its set of vertices,
by~$\Edg(\Gam)$ its set of edges, by~$\Con(\Gam)$ its set of connected
components, and by~$\Hol(\Gam)$ its set of \holes; if no confusion can arise,
we drop the dependence of~$\Gam$.
For a finite set~$X$, let $\FF_2^X$ denote the vector space of
functions $f: X \to \FF_2$ with the values in the $2$-element
field~$\FF_2$.
To any incidence relation $I \subset X \times Y$ between finite
sets~$X$ and~$Y$ we associate a linear map
$\rho = \rho_I: \FF_2^X \to \FF_2^Y$ given by
\begin{equation}
\label{eq:incidence-map}
(\rho f)(y) = \sum_{(x,y) \in I} f(x) \ .
\end{equation}
In particular, there is a natural sequence of linear maps:
\begin{equation}
\label{eq:sequence}
0 \to \FF_2^{\Con} \to \FF_2^{\Ver} \to
\FF_2^{\Edg} \to \FF_2^{\Hol} \to 0 \ ,
\end{equation}
where all the maps are of the form \eqref{eq:incidence-map}
associated to the obvious incidence relations.

\begin{theorem}
\label{thm:CO-MainResult}
Each of the following conditions on a finite graph~$\Gam$ is equivalent
to \COP:
\begin{eqnarray}
\label{eq:edge-numbering}
&&\text{The edges of~$\Gam$ can be linearly ordered so that different}\\
\nonumber
&&\text{\holes in~$\Gam$ will have different maximal edges.}\\
\label{eq:TCO-numerical}
&&|\Hol| = |\Edg| - |\Ver| + |\Con| \ .\\
\label{eq:exact-sequence}
&&\text{The sequence \eqref{eq:sequence} is exact.}
\end{eqnarray}
\end{theorem}

Using~\eqref{eq:exact-sequence}, we obtain the following
corollary implying Proposition~\ref{pr:cyclic-B-to-A}.

\begin{corollary}
\label{cor:product-signs}
In a \CO graph, one can attach the signs $\pm 1$ to all edges in
such a way that the product of these signs along every \hole is
equal to~$-1$.
Furthermore, such an attachment is unique up to simultaneous
changes of signs for all edges incident to a given vertex.
\end{corollary}

\smallskip

As a preparation for the proof of Theorem \ref{thm:CO-MainResult},
we show that certain weaker versions of
\eqref{eq:edge-numbering}-\eqref{eq:exact-sequence} hold for
arbitrary finite graphs.

\begin{proposition}
\label{pr:three-facts}
The following properties hold for an arbitrary finite graph~$\Gam$:
\begin{eqnarray}
\label{eq:edge-numbering-weaker}
&&\text{The edges of~$\Gam$ can be linearly ordered so that an edge}\\
\nonumber
&&\text{is maximal in some \hole if and only if it is one of}\\
\nonumber
&&\text{the last $|\Edg| - |\Ver| + |\Con|$ edges on the list.}\\
\label{eq:TCO-numerical-weaker}
&&|\Hol| \geq |\Edg| - |\Ver| + |\Con| \ .\\
\label{eq:exact-sequence-weaker}
&&\text{The sequence $0 \to \FF_2^{\Con} \to \FF_2^{\Ver} \to
\FF_2^{\Edg} \to \FF_2^{\Hol}$ is exact.}
\end{eqnarray}
\end{proposition}

\proof
To construct a linear ordering of edges satisfying
\eqref{eq:edge-numbering-weaker},
take the first $|\Ver| - |\Con|$ edges (in an arbitrary order)
so that they form a spanning forest in~$\Gam$ (the union of spanning trees
in all connected components).
Then order the remaining $|\Edg| - |\Ver| + |\Con|$
edges so that every new edge connects two vertices at the minimal possible
distance in the graph formed by the preceding edges.
By virtue of this construction, each of the last $|\Edg| - |\Ver| +
|\Con|$ edges becomes the maximal edge in at least one \hole of~$\Gam$.

The bound \eqref{eq:TCO-numerical-weaker} follows at once from
\eqref{eq:edge-numbering-weaker}.

Finally, the only statement in \eqref{eq:exact-sequence-weaker}
which may be not immediately clear is the inclusion
$${\rm Ker}(\FF_2^{\Edg} \to \FF_2^{\Hol}) \subset
{\rm Im}(\FF_2^{\Ver} \to \FF_2^{\Edg}) \ .$$
To prove this inclusion, first consider the case where~$\Gam$ is a forest.
In this case, we need to show that the map $\FF_2^{\Ver} \to
\FF_2^{\Edg}$ is surjective, which is well-known (and follows
easily by induction on the number of edges).
Returning to the general case, choose a linear ordering of the edges
satisfying \eqref{eq:edge-numbering-weaker}.
It is clear that a function $g \in {\rm Ker}(\FF_2^{\Edg} \to
\FF_2^{\Hol})$ is uniquely determined by its restriction to the
first $|\Ver| - |\Con|$ edges.
Since these edges form a forest, there is a function $f \in
\FF_2^{\Ver}$ whose image in $\FF_2^{\Edg}$ agrees with~$g$ on them;
therefore, the image of~$f$ is~$g$, as desired.
\endproof

\smallskip

\noindent{\bf Proof of Theorem \ref{thm:CO-MainResult},
\eqref{eq:edge-numbering}$\Longleftrightarrow$
\eqref{eq:TCO-numerical}$\Longleftrightarrow$
\eqref{eq:exact-sequence}}.
The implication
$\eqref{eq:exact-sequence} \Longrightarrow \eqref{eq:TCO-numerical}$
is clear, while
$\eqref{eq:TCO-numerical} \Longrightarrow \eqref{eq:edge-numbering}$
is immediate from \eqref{eq:edge-numbering-weaker}. Finally, to show
$\eqref{eq:edge-numbering} \Longrightarrow \eqref{eq:exact-sequence}$,
we note that, in view of \eqref{eq:exact-sequence-weaker},
the exactness of the sequence \eqref{eq:sequence}
simply means that its right-most map $\FF_2^{\Edg} \to \FF_2^{\Hol}$
is surjective. The latter condition follows easily
from~\eqref{eq:edge-numbering}, and we are done.
\endproof

To complete the proof of Theorem~\ref{thm:CO-MainResult}, we will
show that the \COP of~$\Gam$ is equivalent to \eqref{eq:edge-numbering}.
First some preparation.
Let~$i$ and~$j$ be two non-adjacent vertices of a finite graph $\Gamma'$,
and let~$\Gam$ be obtained from $\Gamma'$ by adjoining the edge $\{i,j\}$.
By a \emph{chain} connecting~$i$ and~$j$ in~$\Gam'$ we will mean
a sequence of $t \geq 3$ distinct vertices $(i=i_1, i_2, \ldots, i_t=j)$
such that the only edges between them are $\{i_p, i_{p+1}\}$ for
$p = 1, \dots, t-1$.
It is then clear that the \holes in $\Gam$ containing the edge $\{i,j\}$
are in a bijection with the chains connecting~$i$ and~$j$ in~$\Gam'$.
In particular, $i$ and $j$ belong to a \emph{unique} \hole
in~$\Gam$ if and only if they are connected by a \emph{unique}
chain in~$\Gam'$.

\begin{lemma}
\label{lem:lined-up}
Suppose~$\Gam'$ is a \CO graph, and~$i$ and~$j$ are two
non-adjacent vertices of~$\Gam'$ connected by a unique
chain $(i=i_1, i_2, \ldots, i_t=j)$.
Then there exists an orientation of~$\Gam'$
making all \holes cyclically oriented and having the chain
$(i_1, i_2, \ldots, i_t)$ linearly oriented:
$i_1 \to i_2 \to \cdots \to i_t$.
\end{lemma}

\proof
Let~$\Gam''$ be the induced subgraph of~$\Gam'$ on the set
of vertices $\Ver(\Gam') - \{i_1, \dots, i_t\}$.
For $p = 1, \dots, t-1$, let $\Gam_p$ denote the union of all
connected components of~$\Gam''$ joined by an edge with~$i_p$ and~$i_{p+1}$.
We claim that the subgraphs $\Gam_1, \dots, \Gam_{t-1}$ are pairwise
disjoint, and so are disconnected from each other in~$\Gam''$.
Indeed, if $\Gam_p \cap \Gam_q \neq \emptyset$ for some $p < q$
then $i_p$ and $i_{q+1}$ would be connected by a chain having all
interior vertices in~$\Gam''$; but this contradicts the uniqueness of
a chain connecting~$i$ and~$j$.
Clearly, if  an edge $\{i_p, i_{p+1}\}$ belongs to some
\hole~$Z$ in~$\Gam'$, then the rest of the vertices of~$Z$
belong to~$\Gam_p$.

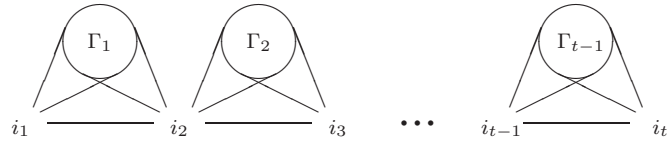
\begin{center}
\begin{picture}(240,55)
\put(0,0){
  \put(0,0){\HBCenter{$\scriptstyle i_1$}}
  \put(60,0){\HBCenter{$\scriptstyle i_2$}}
  \put(10,4){\line(1,0){40}}
  \put(120,0){\HBCenter{$\scriptstyle i_3$}}
  \put(70,4){\line(1,0){40}}
  \multiput(0,0)(60,0){2}{
    \put(4.9,10){\line(2,5){12}}
    \put(55.1,10){\line(-2,5){12}}
    \put(7.5,8.1){\line(2,1){29}}
    \put(52.5,8.1){\line(-2,1){29}}
    \put(30,35){\circle{26}}
    }
  \put(30,35){\HVCenter{$\scriptstyle \Gamma_1$}}
  \put(90,35){\HVCenter{$\scriptstyle \Gamma_2$}}
  \multiput(145,4)(5,0){3}{\circle*{2}}
  \put(180,0){\HBCenter{ $\scriptstyle i_{t-1}$}}
  \put(190,4){\line(1,0){40}}
  \put(240,0){\HBCenter{ $\scriptstyle i_t$}}
  \put(180,0){
    \put(4.9,10){\line(2,5){12}}
    \put(55.1,10){\line(-2,5){12}}
    \put(7.5,8.1){\line(2,1){29}}
    \put(52.5,8.1){\line(-2,1){29}}
    \put(30,35){\circle{28}}
    }
  \put(210,35){\HVCenter{ $\scriptstyle \Gamma_{t-1}$}}
}
\end{picture}
\end{center}

Now take any orientation of~$\Gam'$ making all the \holes
cyclically oriented.
If all the edges $\{i_p, i_{p+1}\}$ are oriented as $i_p \to i_{p+1}$,
there is nothing to prove.
Otherwise, modify the orientation as follows: reverse every directed edge
$i_{p+1} \to i_p$ together with all the edges inside~$\Gam_p$ and
all the edges connecting~$i_p$ and $i_{p+1}$ with~$\Gam_p$.
Clearly, the new orientation still has all the \holes in~$\Gam'$
cyclically oriented, so we are done.
\endproof

\noindent{\bf Proof of Theorem \ref{thm:CO-MainResult},
\eqref{eq:edge-numbering}$\Longrightarrow$ \COP}.
Thus, we assume that~$\Gam$ is a finite graph whose edges are
ordered in accordance to~\eqref{eq:edge-numbering}.
To show that~$\Gam$ is \CO, we proceed by induction on the number of edges~$N$
of~$\Gam$. If~$N=0$, there is nothing to prove; so we assume that~$\Gam$ has at
least one edge.
Let~$e = \{i,j\}$ be the maximal edge in the chosen ordering, and
let~$\Gam'$ be the graph obtained from~$\Gam$ by removing~$e$.
Since~$e$ belongs to at most one \hole in~$\Gam$, it is easy to
see that the \holes in~$\Gam'$ are precisely the \holes in~$\Gam$
that do not contain~$e$.
It follows that~$\Gam'$ satisfies \eqref{eq:edge-numbering}.
By induction, we can assume that~$\Gam'$ is \CO.
Choose an orientation of~$\Gam'$ making all the \holes cyclically oriented.
If~$\Gam$ has no \hole containing~$e$ then $i$ and $j$ belong to
different connected components of~$\Gam'$, and so~$e$ can be
oriented either way, making all \holes in~$\Gam$ cyclically oriented.
This leaves us with the case where~$\Gam$ has a unique \hole
containing~$e$.
Then~$i$ and~$j$ are connected by a unique chain in~$\Gam'$.
By Lemma~\ref{lem:lined-up}, we can assume that this chain is linearly
oriented in a chosen orientation of~$\Gam'$.
Therefore, orienting~$e$ so that the only \hole containing~$e$
becomes cyclically oriented, yields an orientation of~$\Gam$
making all \holes cyclically oriented, thus showing that~$\Gam$ is \CO.
\endproof

\smallskip

It remains to show that a \CO graph satisfies \eqref{eq:edge-numbering}.
Again we need some preparation.
We start with two lemmas.

\begin{lemma}
\label{lem:two-holes}
Let~$\Gam$ be a \CO graph, and suppose two \holes~$Z$
and~$Z'$ share a common edge $\{i,j\}$.
Then~$i$ and~$j$ are the only common vertices of~$Z$ and $Z'$, and
$\{i,j\}$ is the only edge joining a vertex of~$Z$ with a vertex of~$Z'$.
\end{lemma}

\proof
Let~$i, j, i_1, \dots, i_p$ (resp.~$i, j, i'_1, \dots, i'_q$) be all
the vertices of~$Z$ (resp.~$Z'$) written in the cyclical order.
Without loss of generality, we can assume that $i_1 \neq i'_1$
(this can be always achieved by replacing $\{i,j\}$ with another
edge if necessary).

To show that~$i$ and~$j$ are the only common vertices of~$Z$
and~$Z'$, we suppose that $i_k = i'_\ell$ for some $k$ and $\ell$,
and choose the smallest possible~$k$ with this property.
Thus, the sequence of vertices
$j, i_1, \dots, i_k, i'_{\ell - 1}, \dots, i'_1, j$
forms a cycle of length $k~+~\ell$.
Now choose $r \in \{1, \dots, k\}$ and $s \in \{1, \dots, \ell - 1\}$ so that
$i_r$ and $i'_s$ are joined by an edge, and $r+s$ is smallest
possible with this property.
Then the vertices $j, i_1, \dots, i_r, i'_s, i'_{s-1}, \dots, i'_1$
form a \hole, say~$Z''$.
To obtain the desired contradiction, it remains to observe that no
orientation of edges of~$\Gam$ can make all three \holes $Z, Z'$,
and $Z''$ cyclically oriented: indeed, if say~$j$ is a head of $\{i,j\}$
then~$j$ must be a tail of both $\{j,i_1\}$ and $\{j, i'_1\}$, and
so $Z''$ will not be cyclically oriented.
Having proven that the vertices $i_1, \dots, i_k, i'_1, \dots,
i'_\ell$ are all distinct, the same argument as above shows that
there are no edges of the form $\{i_r, i'_s\}$, and
we are done.
\endproof

\begin{lemma}
\label{lem:edge-single-hole-j}
Suppose that in a \CO graph~$\Gamma$, there is at least one edge through a
vertex~$j$.
Then there is an edge through~$j$ that belongs to at most one \hole.
\end{lemma}

\proof
Fix an orientation of~$\Gam$ making all \holes cyclically oriented.
Arguing by contradiction, suppose that every edge through~$j$
belongs to at least two \holes.
Start with an arbitrary edge $\{j, i_1\}$, let $Z_1$ be a \hole
containing this edge, and let $\{j, i_2\}$ be the second edge
through~$j$ that belongs to~$Z_1$.
Then pick a \hole~$Z_2 \neq Z_1$ that contains the edge
$\{j, i_2\}$, and let $\{j, i_3\}$ be the second edge
through~$j$ that belongs to~$Z_2$.
Continuing in the same way, we construct a sequence of \holes
$Z_1, Z_2, \dots$ and a sequence of vertices $i_1, i_2, \dots$
such that each $Z_k$ contains edges $\{j,i_k\}$ and $\{j, i_{k+1}\}$.
Renaming the vertices if necessary, we can assume without loss of
generality that the vertices $i_1, \dots, i_p$ and the \holes
$Z_1, \dots, Z_p$ are pairwise distinct, but
$i_{p+1} = i_1$ and $Z_{p+1} = Z_1$.
In view of Lemma~\ref{lem:two-holes}, we have $p \geq 3$.
Also since all the \holes $Z_1, \dots, Z_p$ must be cyclically
oriented, $p$ is even, and so $p \geq 4$.

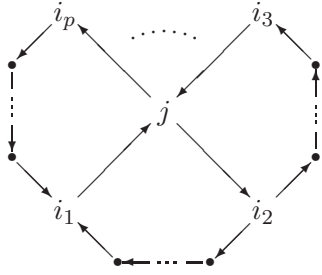
\begin{center}
\begin{picture}(120,110)
\put(60,60){
  \put(0,0){\HVCenter{\small $j$}}
  \put(-37.5,-37.5){\HVCenter{\small $i_1$}}
  \put(-17.5,-57.5){\circle*{3}}
  \put(-19.5,-55.5){\vector(-1,1){13}}
  \put(-5.5,-57.5){\vector(-1,0){10}}
  \multiput(-2.5,-57.5)(3,0){3}{\line(1,0){1}}
  \put(8.5,-57.5){\line(1,0){6}}
  \put(17.5,-57.5){\circle*{3}}
  \put(32.5,-42.7){\vector(-1,-1){13}}
  \put(37.5,-37.5){\HVCenter{\small $i_2$}}
  \put(-32.5,-32.5){\vector(1,1){27.5}}
  \put(5,-5){\vector(1,-1){27.5}}
  \put(57.5,-17.5){\circle*{3}}
  \put(42.5,-32.7){\vector(1,1){13}}
  \put(57.5,-14.5){\line(0,1){6}}
  \multiput(57.5,-5)(0,3){3}{\line(0,1){1}}
  \put(57.5,5.5){\vector(0,1){10}}
  \put(57.5,17.5){\circle*{3}}
  \put(55.5,19.5){\vector(-1,1){13}}
  \put(37.5,37.5){\HVCenter{\small $i_3$}}
  \put(32.5,32.5){\vector(-1,-1){27.5}}
  \put(-57.5,-17.5){\circle*{3}}
  \put(-55.5,-19.7){\vector(1,-1){13}}
  \put(-57.5,-5.5){\vector(0,-1){10}}
  \multiput(-57.5,5)(0,-3){3}{\line(0,-1){1}}
  \put(-57.5,14.5){\line(0,-1){6}}
  \put(-57.5,17.5){\circle*{3}}
  \put(-42.5,32.5){\vector(-1,-1){13}}
  \put(-37.5,37.5){\HVCenter{\small $i_p$}}
  \put(-5,5){\vector(-1,1){27.5}}
  \put(0,30){\circle*{1}}
  \put(4,29.73){\circle*{1}}
  \put(-4,29.73){\circle*{1}}
  \put(7.93,28.93){\circle*{1}}
  \put(-7.93,28.93){\circle*{1}}
  \put(11.72,27.61){\circle*{1}}
  \put(-11.72,27.61){\circle*{1}}
}
\end{picture}
\end{center}

Now let $i_1, \dots, i_2, \dots, i_3, \dots, i_p, \dots, i_1$ be
the sequence of vertices in $Z_1 \cup \cdots \cup Z_p - \{j\}$
written in the natural cyclical order, so that each interval
$i_k, \dots, i_{k+1}$ is a chain $Z_k - \{j\}$.
Choose a shortest possible interval (in the cyclical order) $j_1, \dots, j_q$
in this sequence consisting of $q \geq 3$ pairwise distinct vertices and such
that $\{j_1, j_q\}$ is also an edge. (It is easy to see that such
intervals exist because every three consecutive vertices in our
sequence are distinct.)
The vertices $j_1, \dots, j_q$ obviously form a \hole.
Clearly, this \hole is not a part of any $Z_k - \{j\}$.
Thus, it must contain some $i_k$ as an interior vertex.
But this contradicts the \COP since the two edges through
$i_k$ belong to the \holes $Z_k$ and $Z_{k-1}$ (with the convention
that $Z_0 = Z_p$) and so either both are oriented towards~$i_k$, or
both are oriented away from~$i_k$.
This contradiction proves the result.
\endproof

\noindent{\bf Proof of Theorem \ref{thm:CO-MainResult},
\COP$\Longrightarrow$ \eqref{eq:edge-numbering}}.
Let~$\Gam$ be a cyclically orientable graph
having at least one edge.
By Lemma~\ref{lem:edge-single-hole-j},~$\Gam$ has an edge~$e$
belonging to at most one \hole.
Let~$\Gam'$ be obtained from $\Gam$ by removing the edge~$e$.
While proving the reverse implication, we have already noticed that the \holes in~$\Gam'$ are precisely
the \holes in~$\Gam$ that do not contain~$e$.
It follows that an orientation of~$\Gam$ making all \holes cyclically
oriented, restricts to the orientation of~$\Gam'$ with the same property.
Thus,~$\Gam'$ is \CO.
By induction on the number of edges, there is a linear ordering
of the edges of~$\Gam'$ satisfying \eqref{eq:edge-numbering}.
Adding~$e$ as the maximal element gives a desired linear ordering of
edges for~$\Gam$, finishing the proof of Theorem~\ref{thm:CO-MainResult}.
\endproof

\begin{remark}
\label{rem:checking-4}
Theorem~\ref{thm:CO-MainResult} suggests the following way to
check whether a given skew-symmetrizable matrix satisfies
condition~(4) in Theorem~\ref{thm:MainResult}.
Start by checking whether the edges of (the underlying graph of)
the diagram~$\Gam(B)$ can be ordered in accordance
with~\eqref{eq:edge-numbering}
(we leave aside the question of a practical implementation of such a check).
If \eqref{eq:edge-numbering} cannot be satisfied, then~$B$ is
\emph{not} of finite type.
Otherwise, go over the list of edges ordered in accordance with
\eqref{eq:edge-numbering}, and attach the signs to them in the
following way.
If a current edge~$\{i,j\}$ is not maximal in any \hole, choose
the sign~$\varepsilon_{ij}$ arbitrarily; otherwise,~$\{i,j\}$ is
maximal in a unique \hole~$Z$, and we determine~$\varepsilon_{ij}$
from the condition that the product of signs along the edges
of~$Z$ is equal to~$-1$.
Finally, define a \QC companion~$A$ of~$B$ by setting
$A_{ij} = -\varepsilon_{ij} |B_{ij}|$ for all~$i \neq j$.
Then~$B$ satisfies~(4) and so is of finite type if and only if~$A$
is positive (the latter condition can be checked for example by
using the Sylvester criterion).
\end{remark}

\section*{Acknowledgments}
The first solution of the recognition problem for cluster
algebras of finite type was given by A.~Seven in \cite{Seven}
as a part of his Ph.D. thesis written under the supervision of
A.~Zelevinsky at Northeastern University.
M.~Barot and C.~Geiss started working on an alternative solution after
attending A.~Seven's talk at the Sixth International Joint AMS-SMM Meeting
in May 2004.
They joined forces with A.~Zelevinsky after he attended
their presentation of a preliminary version of Theorem~\ref{thm:MainResult}
at the XI International Conference on Representations of Algebras, Mexico,
August 2004.

The authors thank Alex Postnikov for stimulating discussions and
help in proving the results in Section~\ref{sec:cog}.


\begin{thebibliography}{xxx}

\bibitem{B}{N.~Bourbaki, \sl{Groupes and alg\`ebres de Lie}, Ch. 4,5
    et 6. Hermann, Paris 1968.}

\bibitem{FZ1}{S.~Fomin, A.~Zelevinsky, \sl{Cluster algebras I:
      Foundations}. J. Amer. Math. Soc. {\bf 15} (2002), no. 2, 497-529.}

\bibitem{FZ2}{S.~Fomin, A.~Zelevinsky, \sl{Cluster algebras II: Finite type
    classification}. Invent. Math. {\bf 154} (2003), no 1. 63-121.}

\bibitem{kac}{V.~Kac, {\sl Infinite dimensional Lie algebras}, 3rd edition,
Cambridge University Press, 1990.}

\bibitem{Ri}{C.~M.~Ringel, {\sl Tame algebras and integral quadratic forms,}
Lecture Notes in Mathematics, {\bf 1099},
Springer-Verlag, Berlin, 1984.}

\bibitem{Seven}
{A.~Seven, \sl{Recognizing cluster algebras of finite type,}
\texttt{math.CO/0406545}.}
\end{thebibliography}
\end{document}